\documentclass[11pt]{article}
\usepackage{amsfonts}
\addtocounter{page}{0}
\usepackage[german]{babel}
\begin{document}

\begin{center}
\large{\bf Another look at {\it e}\footnote{The initial version of this paper was
submitted for publication on July12, 2009}}

\vskip .3in

Samuel L. Marateck\footnote{{\it email: marateck@courant.nyu.edu}}\\ 

{\it Courant Institute of Mathematical Sciences, New York University, New York,
  N.Y. 10012}
\end{center}

\vskip .3in
\footnoterule
{\bf Abstract}
\begin{quotation}
\noindent
\small{This note describes  a way of obtaining {\it e} that differs from the
  standard one. It could be used as an alternate way of showing how the value
of {\it e} is obtained. No attempt is made to show the existence of the limit
in the definition  of {\it e} that appears in the final equation.
}
\end{quotation}

\footnoterule
\vskip .3in

\noindent
{\bf 1. Introduction.} \\
\noindent
Traditionally the value of $e$ has been obtained, for instance, by taking the
limit of ever-decreasing interest intervals in the compound interest formula
(see Greenleaf [1]) or linear interpolation (see Flanders and Price [2]). We
describe an alternate technique of obtaining $e$ that should have pedagogic
value. In this section we give an approximation of $e$ using this technique and
generalize it in the next section.

\noindent
If $f'(x)$, the derivative of $f(x),$ exists at point $x$, and you start at
point $x$ and move a distance $\Delta x$, the value at the point $x + \Delta$
is given by

\begin{equation} f(x + \Delta x) \cong f(x) + f'(x)\cdot \Delta x \end{equation}

\noindent
We want to find a constant, let's call it $e$, such that when it's raised to
the power $x$ obtaining the function $e^x$, the function's derivative is also
$e^x$\footnote{Our analysis also holds if $f(x) = Ce^x$ where $C$ is a
constant.}.

\noindent
Since $f'(x)$ equals $f(x),$ we rewrite equation (1) as

\begin{equation} f(x + \Delta x) \cong f(x)(1 + \Delta x) \end{equation}

\noindent
We will analyse this in the interval [1,2].  Let's take x = 1 and $\Delta x$ =
0.1. So $x + \Delta x$ is 1.1. Equation (2) gives

      \begin{equation}f(1.1) \cong f(1)(1 + 0.1) \end{equation}

\noindent
or

      \begin{equation}e^{1.1} \cong 1.1e \end{equation}.

\noindent
Now take $x = 1.1$ and use the same value of $\Delta x$, i.e., 0.1. 
We will be using the same increment in $x$ in this and all subsequent steps
since eventually we will let $\Delta x$ approach zero. Continuing in this way

     \begin{equation}f(1.1 + 0.1) \cong 1.1f(1.1) \end{equation}

\noindent
So $f(1.2) \cong 1.1e^{1.1}$. Or

   \begin{equation}e^{1.2} \cong (1.1)^2e \end{equation}

\noindent
Eventually we will get $e^2$ on the left side of the equation, so we can solve
for e.  So let's compute $e^{1.3}$. We get $e^{1.3} \cong 1.1e^{1.2}$ But this
equals $(1.1)^3e$.  If we extrapolate to $x = 1.8$, we see that

     \begin{equation}e^{1.9} \cong (1.1)^9e \end{equation}

\noindent
and finally that

      \begin{equation}e^2 \cong (1.1)^{10}e \end{equation}

\noindent
Solving for $e$ we get $e \cong (1.1)^{10}$ or $e$ equals 2.59 to three digits,
where the 10 corresponds to dividing 1 by 0.1.  Equation (1) presupposes that
$\Delta x$ approaches zero.  If we let $\Delta x= .00000001$,
or $10^{-8}$, we raise $(1 + .00000001)$ to $10^{8}$. The answer for $e$ is
2.71828 to five significan figures. 

\noindent
{\bf 2. Generalization.} \\

\noindent
We now sketch the steps that describe the preceding method in general. Using
equation (2), and setting $x = 1$, we write

\begin{equation}e^{1 + \Delta x} \cong e(1 + \Delta x) \end{equation}

\noindent
We continue, letting $x = x+ \Delta x$ and keeping $\Delta x$ the same, and
write

\begin{equation}e^{1 + \Delta x + \Delta x} \cong e^{1 + \Delta x}(1 + \Delta
  x) \end{equation}

\noindent
or

\begin{equation}e^{1 + 2\Delta x} \cong e(1 + \Delta x)^2 \end{equation}

\noindent
We have to add $\Delta x$ to $x$ $1/\Delta x$ times to get $e^2$ on the left
side of these equations. So we get

\begin{equation}e^{1 + (1/ \Delta x) \cdot \Delta x} \cong e(1 + \Delta x)^{1/
\Delta x} \end{equation}

\noindent
or

\begin{equation}e^2 = e(1 + \Delta x)^{1/\Delta x} \end{equation}

\noindent
Solve for e and since the definition of the derivative in equation (1)
lets $\Delta x \to 0$, take the same limit here. We get

\begin{equation}e = \displaystyle \lim_{\Delta x \to 0}(1 + \Delta x)^{1/
\Delta x} \end{equation}

\noindent
which is one of the definitions of $e$.

\noindent
{\bf 3. References.} \\

\noindent
[1] Greenleaf, Fredrick E.. {\it Quantatative Reasoning, 3/e}, McGraw Hill
(2006).

\noindent
[2] Flanders, Harley and Price, Justin J., {\it Calculus with Analytic
  Geometry}, Academic Press (1978).

\end{document}